\documentclass[oneside,leqno,10pt]{article}

\newif\ifpdf
\ifx\pdfoutput\undefined
	\pdffalse	
\else
	\pdfoutput=1	
	\pdftrue
\fi

\usepackage{a4}

\usepackage{oldgerm}
\usepackage{bbold}
\usepackage{amssymb,amsmath,amscd}
\ifpdf
	\usepackage[pdftex]{graphicx}
\else
	\usepackage{graphicx}
\fi
\usepackage[latin1]{inputenc}
\usepackage[ansi]{umlaute}
\usepackage{fancyvrb}
\usepackage{fancyhdr}
\usepackage{pslatex}
\usepackage{makeidx}
\usepackage{epic}
\usepackage{eepic}
\usepackage{nicefrac}

\newcommand{\transpose}[1]{\,\mbox{}^t\hspace{-2pt}{#1}}
\newcommand{\sqmatrixtwo}[4]{\begin{pmatrix} #1 & #2 \\ #3 & #4 \end{pmatrix}}
\newcommand{\smallsqmatrixtwo}[4]{\bigl(\begin{smallmatrix}#1&#2 \\ #3&#4\end{smallmatrix}\bigr)}

\newcommand{\suchthat}{\mbox{\ \ :\ \ }}
\newcommand{\where}{\big|}

\newcommand{\refTits}[1]{\cite[\ref{#1}]{Tits}}

\DeclareMathAlphabet{\mathobb}{U}{bbold}{m}{n}
\newcommand{\UnitMx}{\mathobb{1}}

\newcommand{\diag}{\ensuremath{\operatorname{diag}}}
\newcommand{\tr}{\ensuremath{\operatorname{tr}}}
\newcommand{\id}{\ensuremath{\operatorname{id}}}

\newcommand{\N}{{\mathbb{N}}}
\newcommand{\Z}{{\mathbb{Z}}}
\newcommand{\Q}{{\mathbb{Q}}}
\newcommand{\R}{{\mathbb{R}}}
\newcommand{\C}{{\mathbb{C}}}

\newcommand{\GL}{\operatorname{GL}}
\newcommand{\SL}{\operatorname{SL}}
\newcommand{\Sp}{\operatorname{Sp}}
\newcommand{\Sym}{\operatorname{Sym}}
\newcommand{\Stab}{\operatorname{Stab}}

\newcommand{\D}{\textfrak{D}}
\newcommand{\cP}{\ensuremath{{\mathcal P}}}

\renewcommand{\S}{\mathfrak{S}}
\newcommand{\PolMx}{\Delta}
\newcommand{\dsum}[2]{\ensuremath{d_{#1:#2}}}
\newcommand{\Dsum}[2]{\ensuremath{D_{#1:#2}}}
\newcommand{\esum}[2]{\ensuremath{e_{#1:#2}}}
\newcommand{\ssum}[2]{\ensuremath{s_{#1:#2}}}

\newcommand{\bsum}[2]{\ensuremath{b_{#1:#2}}}
\newcommand{\Bsum}[2]{\ensuremath{B_{#1:#2}}}
\newcommand{\csum}[2]{\ensuremath{c_{#1:#2}}}
\renewcommand{\L}{{\mathfrak{L}}}
\newcommand{\LL}{{\mathbb{L}}}

\newcommand{\A}{{\mathcal A}}
\newcommand{\pol}{_{\text{pol}}}
\newcommand{\lev}{^{\text{lev}}}
\newcommand{\Apol}{\ensuremath{\A\pol}}
\newcommand{\Apollev}{\ensuremath{\A\pol\lev}}
\newcommand{\Gampol}{\ensuremath{\tilde\Gamma\pol}}
\newcommand{\Gampollev}{\ensuremath{\tilde\Gamma\pol\lev}}
\newcommand{\GampolConj}{\ensuremath{\Gamma\pol}}
\newcommand{\GampollevConj}{\ensuremath{\Gamma\pol\lev}}

\newcommand{\Gp}[1]{\ensuremath{\tilde\Gamma_{\text{pol,$#1$}}}}
\newcommand{\Gpl}[1]{\ensuremath{\tilde\Gamma_{\text{pol,$#1$}}^{\text{lev}}}}
\newcommand{\GpC}[1]{\ensuremath{\Gamma_{\text{pol,$#1$}}}}
\newcommand{\GplC}[1]{\ensuremath{\Gamma_{\text{pol,$#1$}}^{\text{lev}}}}

\newcommand{\GamG}{\ensuremath{\Gamma_G}}
\newcommand{\Cpollev}{\ensuremath{{\mathcal C}\pol\lev}}

\newcommand{\ord}{\operatorname{ord}}
\newcommand{\chisym}{\ensuremath\chi^{\text{sym}}}

\newcommand{\DP}{\ensuremath{\mathbb{D}(\PolMx)}}
\newcommand{\Dd}{\ensuremath{\mathbb{D}(\PolMx_d)}}
\newcommand{\De}{\ensuremath{\mathbb{D}(\PolMx_e)}}

\newcommand{\PdashQ}{\ensuremath{P'_{\Q}}}
\newcommand{\Pdash}[2]{\ensuremath{P'_{#1}(#2)}}

\ifpdf
\fi

\newcommand{\mathtype}{undefMathType}
\newcounter{mathcount}[section]
\renewcommand{\themathcount}{\mathtype~\arabic{section}.\arabic{mathcount}}

\newenvironment{bew}[1]{\smallskip\noindent{\bf Proof#1.\\}\rm\renewcommand{\mathtype}{Proof}}
  {\mbox{}\hfill$\square$\medskip} 

\newenvironment{mathenv}[1]{\refstepcounter{mathcount}\bigskip\par\noindent{\bf \themathcount#1.\\}}
  {\par} 

\newenvironment{mathenv*}[1]{\refstepcounter{mathcount}\bigskip\par\noindent{\bf \themathcount#1.}}
  {} 

\newenvironment{defi}[1]{\renewcommand{\mathtype}{Definition}\begin{mathenv}{#1}}
  {\end{mathenv}\medskip}
\newenvironment{defi*}[1]{\renewcommand{\mathtype}{Definition}\begin{mathenv*}{#1}}
  {\end{mathenv*}}

\newenvironment{lem}[1]{\renewcommand{\mathtype}{Lemma}\begin{mathenv}{#1}}
  {\end{mathenv}}
\newenvironment{lem*}[1]{\renewcommand{\mathtype}{Lemma}\begin{mathenv*}{#1}}
  {\end{mathenv*}}

\newenvironment{kor}[1]{\renewcommand{\mathtype}{Corollary}\begin{mathenv}{#1}}
  {\end{mathenv}}
\newenvironment{kor*}[1]{\renewcommand{\mathtype}{Corollary}\begin{mathenv*}{#1}}
  {\end{mathenv*}}

\newenvironment{satz}[1]{\renewcommand{\mathtype}{Theorem}\begin{mathenv}{#1}}
  {\end{mathenv}}
\newenvironment{satz*}[1]{\renewcommand{\mathtype}{Theorem}\begin{mathenv*}{#1}}
  {\end{mathenv*}}

\newenvironment{anm}[1]{\renewcommand{\mathtype}{Remark}\begin{mathenv}{#1}}
  {\end{mathenv}}
\newenvironment{anm*}[1]{\renewcommand{\mathtype}{Remark}\begin{mathenv*}{#1}}
  {\end{mathenv*}}

\newenvironment{beh*}[1]{\renewcommand{\mathtype}{Proposition}\begin{mathenv*}{#1}}
  {\end{mathenv*}}

\newenvironment{notn}{\renewcommand{\mathtype}{Notation}\begin{mathenv}{}}
  {\end{mathenv}}
\newenvironment{notn*}{\renewcommand{\mathtype}{Notation}\begin{mathenv*}{}}
  {\end{mathenv*}}

\newenvironment{bsp}{\renewcommand{\mathtype}{Example}\begin{mathenv}{}}
  {\end{mathenv}}

\begin{document}

\begin{center}
\Huge {\sffamily\bfseries The Kodaira dimension of\\Siegel modular varieties of genus 3 or higher}\\
Eric Schellhammer\\[.5em]
\end{center}

\section*{Abstract}

We consider the moduli space $\Apol(n)$ of (non-principally) polarised abelian varieties of genus $g\geq3$ with coprime polarisation and full level-$n$ structure.
Based upon the analysis of the Tits building in \cite{Tits}, we give an explicit lower bound on $n$ that is sufficient for the compactified moduli space to be of general type if one further explicit condition is satisfied.

\section{Introduction}

For positive integers $d_1,\dots,d_{g-1}$ define $\dsum{i}{j}:=\prod_{k=i}^jd_k$ where the empty product equals 1. A polarisation of type $(1,d_1,\dots,\dsum{1}{g-1})$ is called coprime if $\gcd(d_i,d_j)=1$ for all $i\neq j$.
Let $\PolMx:=\diag(1,d_1,\dsum{1}{2},\dots,\dsum{1}{g-1})$, $\Lambda:=\smallsqmatrixtwo{0}{\PolMx}{-\PolMx}{0}$. Let $\L:=\Z^{2g}\subset\C^g$ and denote the lattice dual to $\L$ with respect to the bilinear form $\langle x,y\rangle:=x\Lambda\transpose{y}$ by $\L^\vee$.

Recall that the paramodular groups without, with canonical or with full level structure can be defined as follows:
\begin{align*}
\GpC{d} & :=\{M\in\SL(2g,\Q) \where M\Lambda\transpose{M}=\Lambda\}\\
\GplC{d} & :=\{M\in\GpC{d} \where M|_{\L^\vee/\L}=\id|_{\L^\vee/\L}\}\quad\text{and}\\
\GpC{d}(n) & :=\{M\in\GpC{d} \where M\equiv\UnitMx\mod n\}.
\end{align*}
When it is obvious which polarisation we refer to we simply write $\GampolConj, \GampollevConj$ and $\GampolConj(n)$, respectively. All these groups act on the Siegel upper half space $\S_g$ by
$$\smallsqmatrixtwo{A}{B}{C}{D}:\tau\mapsto(A\tau+B)(C\tau+D)^{-1}.$$
The quotient spaces $\Apol$, $\Apollev$ and $\Apol(n)$ are the moduli spaces of Abelian varieties with fixed polarisation of the given type without, with canonical or with full level structure, respectively.

The Kodaira dimension of these spaces is defined via their compactifications. To obtain these we use the method of toroidal compactification introduced in \cite{AMRT}. Several of these spaces have been thoroughly investigated. For principal polarisations, the work of Freitag, Igusa, Mori, Mukai, Mumford, Tai and a number of other authors gives an almost complete picture which of these spaces are rational, unirational or of general type. The only space where the Kodaira dimension could not yet be determined is $\A_6$, i.\,e.~the case $g=6$ without level structure.

For $g=2$ and a polarisation of type $(1,p)$ we know that $\Apol$ is of general type for all primes $p\geq73$ by recent work of Sankaran and Erdenberger. Several other results are known for polarisations of type $(1,t)$ with small $t$. Furthermore, Hulek showed in \cite{H} that $\A_{(1,t)}(n)$ is of general type for $n\geq4$ when $\gcd(n,t)=1$.

However, the analysis of moduli spaces of genus 3 or higher appears to be more complicated.
Tai showed that for $g\geq16$ all these spaces are of general type, but not much is known for lower $g$.
In this paper we want to give a result concerning the cases $g\geq3$, which for $g\geq16$ is weaker than the result by Tai but closes the gap $3\leq g\leq15$.

\begin{center}
{\bf Theorem:} Let $g\geq3$ and let $(1,d_1,\dsum{1}{2},\dots,\dsum{1}{g-1})$ be a coprime polarisation with $\dsum{1}{g-1}\neq2$. Then $\Apol(n)$ is of general type, provided $\gcd(n,\dsum{1}{g-1})=1$, $n\geq3$ and
$$n>\frac{(2^g+1)\dsum{2}{g-2}}{(g+1)2^{g-3}}\min\{C'(\LL_1),C'(\LL_2)\}$$
where $\LL_1=(d_1,\dots,d_{g-1})$, $\LL_2=(d_{g-1},\dots,d_1)$ and
$$C'(\LL(x_1,\dots,x_{g-1})=x_1\max\left\{1,\frac{1}{\sqrt{3}}\max_{2\leq r\leq g}\left\{\sqrt[r]{\prod_{i=1}^{r-1}x_i^i}\right\}\right\}.$$
\end{center}

We prove the appropriate behaviour of $h^0(K^k)$ by relating it to the line bundle $L$ of modular forms of weight~1 and then using Hirzebruch proportionality.
For this relation we need a cusp form with respect to $\GampolConj(n)$, which we denote by $\overline\chi$ and construct from the $\Sp(2g,\Z)$-cusp form given by the product of all even theta values.
The space $\Apol(n)$ is given a toroidal compactification and on this $\overline\chi$ can be extended to the boundary.
The weight of $\overline\chi$ and its order of vanishing on the boundary are calculated by analysing the maps used to construct it.

When describing toroidal compactification, our notation is based on \cite[Section~3C]{HKW}. In particular, $\cP(F)\subset\Sp(2g,\R)$ is the stabiliser of a rational boundary component $F$, $\cP'(F)$ is the centre of the unipotent radical of $\cP(F)$ and $\Pdash{\Gamma}{F}:=\cP'(F)\cap\Gamma$ its relevant lattice part, where $\Gamma$ is any of the above groups.

\label{sectBnC}

\section{Vanishing on the boundary}

In a toroidal compactification the boundary is composed of several different parts which correspond to rational polyhedral cones in the closure $\overline C$ of the cone of positive definite, symmetric $g\times g$ matrices. The (open) boundary components of codimension~1 correspond\footnote{This is not a 1-to-1-correspondence, since we have to consider several copies of $\overline C$.} to 1-dimensional cones (i.\,e.~rays) in $\overline C$. If the ray is generated by a matrix of rank 1 (which implies that it lies on the boundary $\overline C\backslash C$) we call the corresponding rational boundary component a \emph{corank-1 boundary component}.

These corank-1 boundary components play a crucial part in determining the order of vanishing of a cusp form on all of the boundary.
In the principally polarised case this is shown using the result by Barnes and Cohn in \cite{BC}.

For the non-principally polarised case this theorem unfortunately cannot be established; in fact, there is a counterexample which we will give in \ref{counterBnC}. Nevertheless, a generalisation of the result
by Barnes and Cohn provides a weaker bound which may be used instead.

Following the paper \cite{BC} we generalise their theorem 3 to
some more general lattices which correspond to the non-principally
polarised case. We first recall some notation.

\begin{notn}{}
Let $f(x):=xA\transpose{x}$ and $h(x):=xB\transpose{x}$ be two quadratic
forms with real symmetric matrices $A$ and $B$, and define their inner
product as $(f,h):=\tr(AB):=\sum_{i,j}a_{ij}b_{ij}.$
For positive definite $f$ denote by $M(f)$ its arithmetic minimum,
i.\,e.\ the minimum of $f(x)$ with integral $x\neq0$.
If $\LL$ is a lattice of matrices, we shall write $f\in\LL$ to denote that $f$ can be given as above with $A\in\LL$.
\end{notn}

The theorem by Barnes and Cohn is used in the context of moduli of principally polarised abelian varieties in form of the following
\begin{satz}{}
\label{BnCused}
Let $f$ be a real positive definite $n$-ary form and denote by $\LL^0$ the lattice of all positive definite or positive semi-definite integral forms and by $\LL_1\subset \LL^0$ the sublattice of forms of rank 1. Then 
$$\min_{h\in \LL^0\backslash\{0\}}(f,h)\geq\min_{h\in \LL_1}(f,h).$$
Furthermore, there always exists a form of rank 1 realising this minimum.
\end{satz}
\begin{bew}{}
This is an immediate consequence of \cite[Theorem~3]{BC}.
\end{bew}

The main connection between extending pluricanonical forms to a toroidal compactification and this corollary is \cite[Chapter~IV, paragraph~1, Theorem~1]{AMRT}.
The precise correspondence between the vanishing on the corank-1 boundary components and on the rest of the boundary is\footnote{See also \cite[Theorem~1.1]{Tai}.}:

\begin{kor*}{}\label{Tai1.1}\\
Suppose $\D=\S_g$ and $\Gamma=\Sp(2g,\Z)$. Let $\chi$ be an automorphic form of weight $l(g+1)$ with respect to $\Gamma$, $\omega=\bigwedge_{i\leq j}d\tau_{i,j}$, $\chi\omega^{\otimes l}\in\Omega^N(\S_g/\Gamma)^{\otimes l}$, and let $\overline{\S_g/\Gamma}^0$ be the smooth part of a toroidal compactification of $\S_g/\Gamma$.
Then
$$\chi\omega^{\otimes l}\text{extends to }\overline{\S_g/\Gamma}^0
\iff
\left\{\begin{gathered}
\text{$\chi$ vanishes on all}\\\text{rational corank-1 boundary components}\\\text{of order at least $l$.}
\end{gathered}\right.$$
\end{kor*}

\begin{bew}{}
The proof of this statement can be found in \cite{AMRT}. But although the theorem as stated here is only valid for the principally polarised case, the proof for the non-principally polarised case differs from this one only in the substitution of \ref{BnCused} by a generalisation. Therefore, we want to sketch the proof in order to show how the reduction to forms of rank~1 can be achieved.

Since we have a principal polarisation, $(P'(F))^\vee$ consists of integer matrices for all rational boundary components $F$. Therefore, according to \ref{BnCused}, the minimum of $(f,h)$ with $f\in(P'(F))^\vee$ over all $h\in P'(F)\cap\overline{C(F)}$ is obtained for a form $h$ of rank 1, where $C(F)$ is the self-adjoint cone corresponding to $F$. For any such $h$ we can find a corank-1 boundary component $F_1\prec F$ with $h\in P'(F_1)\cap\overline{C(F_1)}$. Because the coefficients $a_f^F$ of the Fourier-Jacobi expansion are the same for every pair $F\succ F_1$ we can now bound the minimum over all $h$ for all $F$ by the minimal order of vanishing on all rational corank-1 boundary components.
\end{bew}

\subsection{Non-principal polarisations}

\ref{BnCused} depends heavily on the fact that we consider the minimum over {\em all} integral forms $h$. However, this is only the case if we apply it to principal polarisations. Otherwise the matrix of the bilinear form
$h$ is no longer simply an element of $\Sym(g,\Z)$ but of a sublattice.
To make things precise we define the relevant lattices as follows.

\begin{defi}{: Tits Lattice}
In \cite[Paragraph~3D]{HKW} a standard rational boundary component corresponding to the lattices of rank $g$ is defined and denoted by $F^{(0)}$. (This is yet independent of $\Gamma$.)
By the {\em Tits lattice} we mean the lattice $\LL=P'(F^{(0)})\cap\overline{C(F^{(0)})}$ where we identify the containing space $\cP'(F^{(0)})$ with the space of symmetric matrices as in \cite[Paragraph~3D]{HKW}. If the type of the polarisation is given by $(1,d_1,\dots,\dsum{1}{g-1})$ and we have no level structure we also write $\LL(1,d_1,\dots,\dsum{1}{g-1}).$
\end{defi}

\begin{anm}{}
The definition of the Tits Lattice only considers the standard corank-$g$ boundary component. However, \refTits{onlyoneh} tells us that this is no restriction since for square-free, coprime polarisations all corank-$g$ boundary components are conjugate under the action of $\GampolConj$.
\end{anm}

\begin{defi*}{: Characteristic values of a lattice}
\begin{itemize}
\item Let $\LL\subset\Sym(n,\Z)$ be a sublattice of the lattice of symmetric matrices
and define the subsets
$\LL^0\subset \LL$ and $\LL^+\subset \LL^0$
of positive semi-definite (including the zero matrix) and positive definite matrices, respectively.
Let $\LL_1\subset \LL$ be the subset of rank 1 matrices.
\item If $\LL$ is of maximal rank, define two characteristic values for the lattice, namely the greatest common divisor of all (non-zero) determinants
\begin{align*}
\mu(\LL) & := \max\{\lambda\in\N\where\forall B\in \LL^+\suchthat\lambda|\det(B)\}
\intertext{and the least value $\nu$ that makes sure that all matrices $\nu C$ are members of the lattice}
\nu(\LL) & := \min\{\lambda\in\N\where\forall C\in\Sym(n,\Z), C\mbox{ positive semi-definite}\suchthat \lambda C\in \LL^0\}.
\end{align*}
\end{itemize}
\end{defi*}

\begin{lem}{}
\label{latticetopol}
The Tits lattice of a polarisation of type $(1,d_1,\dots,\dsum{1}{n-1})$ without level-structure is
$$\LL(1,\dots,\dsum{1}{g-1})=\Big\{M\in\begin{pmatrix}\Z&d_1\Z&\dots&\dsum{1}{n-1}\Z\\d_1\Z&d_1\Z& & \dsum{1}{n-1}\Z\\\vdots & & \ddots & \vdots \\ \dsum{1}{n-1}\Z & \dsum{1}{n-1}\Z & \dots & \dsum{1}{n-1}\Z \end{pmatrix}\where M\mbox{ symmetric}\Big\}$$
and it has the characteristics $\mu(\LL)=\prod_i d_i^{n-i}$ and $\nu(\LL)=\dsum{1}{n-1}.$
\end{lem}
\begin{bew}{}
From \cite[Paragraph~3D]{HKW} we know $$\cP'(F^{(0)})\simeq\{\smallsqmatrixtwo{\UnitMx}{S}{}{\UnitMx}\where S\in\Sym(g,\R)\}\simeq\Sym(g,\R)$$
for the standard rational boundary component $F^{(0)}$.
This isomorphism maps a matrix $M\in\cP'(F^{(0)})$ onto its upper right quarter.
Since $P'(F^{(0)})=\cP'(F^{(0)})\cap\GampolConj$ we are only interested in the symmetric $g\times g$ matrices satisfying the conditions on the upper right quarter of the matrices in $\GampolConj.$ \refTits{CongGampolConj} gives the condition claimed.
\end{bew}

Now we want to give the aforementioned counterexample to the inequality in \ref{BnCused}:
\begin{bsp}
\label{counterBnC}
Let $\LL=\LL(1,17)$ and
$$f(x)=x\sqmatrixtwo{3}{-\frac{14}{17}}{-\frac{14}{17}}{\frac{4}{17}}\transpose{x}\in \LL^\vee.$$
We claim that $\min_{h\in \LL_1}(f,h)=3.$ To show this, define $h_0$ to be a rank 1 form realizing the minimum and let the form be given by the matrix $\smallsqmatrixtwo{a^2}{ab}{ab}{b^2}$. For $h_0\in \LL_1$ we need $17|ab$ and $17|b^2$.
Since the rank of $h_0$ is 1, we cannot have $a=b=0$. If $a=0$ or $b=0$ we obtain
$$(f,h_0)=\tr(fh_0)=\tfrac{4}{17}b^2=4\quad\text{or}\quad(f,h_0)=\tr(fh_0)=3a^2=3,$$
respectively, since 17 divides $b^2$ and the minimality of $h_0$. Hence, $\min_{h\in \LL_1}(f,h)\leq 3.$

Now assume that $ab\neq0$ and $\tr(fh_0)<3$. Since $h_0$ is positive semi-definite, we have $a^2,b^2\in\N$ and hence $a,b\in\R$. Fix $a\in\R$ and define
$$f_a(b):=\tr(fh_0)=3a^2-\tfrac{28}{17}ab+\tfrac{4}{17}b^2=\tfrac{4}{17}(b-\tfrac{7}{2}a)^2+\tfrac{2}{17}a^2.$$
Then $f_a$ has no zeroes and assumes its minimum over $\R$ at $b=\frac{7}{2}a$. Since the assumption that $f_a(\frac{7}{2}a)=\frac{2}{17}a^2<3$ leads to $a^2<\frac{51}{2}$ and we have seen that $a^2\in\N$, this leaves only 10 possible values for $a$.

If $a=\pm1,\pm2$ then the condition $ab\in17\Z$ leads to $b=17b'$ with $b'\in\Z$. Easy calculation shows that $f_a(17b')=3a^2-28ab'+68b'^2\geq3.$ If $a=\pm\sqrt{2},\pm\sqrt{3},\pm\sqrt{5}$ the condition $ab\in17\Z$ leads to $b=17ab'$ with $b'\in\Z$. But now $f_a(17ab')=a^2f_1(17b')\geq3a^2>3.$ Hence, $min_{h\in \LL_1}(f,h)\geq3$, which shows the claim.

On the other hand, for the rank 2 form $h$ with matrix $\smallsqmatrixtwo{6}{17}{17}{51}$ we calculate $(f,h)=\tr(fh)=2,$ so obviously $$\min_{h\in \LL}(f,h)\leq 2<3=\min_{h\in \LL_1}(f,h)$$
which shows that the inequality of \ref{BnCused} cannot be established for the non-principally polarised case with $p=17$.
\end{bsp}

\subsection{Barnes and Cohn generalised}

A generalisation can be achieved if one allows a factor in the inequality which depends on the characteristic values of the given lattice as follows:
\begin{satz}{}
\label{PosDefThm}
Let $f$ be a real positive definite $n$-ary form where $n\geq2$. Then
$$\min_{h\in \LL^+}(f,h)\geq\frac{\sqrt{3}\sqrt[n]{\mu(\LL)}}{\nu(\LL)}\min_{h\in \LL_1}(f,h).$$
\end{satz}
\begin{bew}{}
It is easy to derive $(f,h)\geq\frac{n}{\gamma_n}\sqrt[n]{\mu(\LL)}M(f)$ for all positive definite forms $f,h$ with $h\in \LL^+$ in the manner of \cite[Corollary~1]{BC}.
We obtain $(f,h)\geq\sqrt{3}\sqrt[n]{\mu(\LL)}M(f)$ in the same way as in \cite[Theorem~2]{BC}.
Now, we chose $h_0$ of rank 1 such that $(f,h_0)=M(f)$ and obtain
\begin{align*}
\min_{h\in \LL^+}(f,h) & > \sqrt{3}\sqrt[n]{\mu(\LL)}M(f)
  = \sqrt{3}\sqrt[n]{\mu(\LL)}(f,h_0) \\
 & = \frac{\sqrt{3}\sqrt[n]{\mu(\LL)}}{\nu(\LL)}(f,\nu(\LL)h_0)
  \geq \frac{\sqrt{3}\sqrt[n]{\mu(\LL)}}{\nu(\LL)}\min_{h\in \LL_1}(f,h)
\end{align*}
since $\nu(\LL)h_0\in \LL_1$ from the definition of $\nu(\LL)$ and the fact that $h_0$ has rank 1.
\end{bew}

\begin{kor}{: $(1,t)$-polarisation}
Let $t\in\N, t\geq3$ and $\LL=\LL(1,t)$. Then
$$\min_{h\in \LL^0\backslash\{0\}}(f,h)\geq\sqrt{\frac{3}{t}}\min_{h\in \LL_1}(f,h).$$
\end{kor}
\begin{bew}{}
This follows from \ref{PosDefThm} using the values given in \ref{latticetopol}.
\end{bew}

Unfortunately, for a general lattice of higher dimension
it is not as easily possible to compare the two minima. Nevertheless, for the special case of Tits lattices we can obtain the following theorem:

\begin{satz*}{}\label{BnCgeneral}\\
Let $f$ be a real positive $n$-ary form with $n\geq2$ and let $\LL=\LL(1,d_1,\dots,\dsum{1}{n-1})$. Then
\begin{gather}
\min_{h\in \LL^+}(f,h) \geq \frac{\sqrt{3}}{\sqrt[n]{\prod_{i=1}^{n-1} d_i^i}}\min_{h\in \LL_1}(f,h)\quad\text{and}\label{BnCpart1}\\
\min_{h\in \LL^0\backslash\{0\}}(f,h) \geq C(\LL)\min_{h\in \LL_1}(f,h)\quad\text{where}\label{BnCpart2}\\
 C(\LL) := \min\Big\{1,\min_{2\leq r\leq n}\frac{\sqrt{3}}{\sqrt[r]{\prod_{i=1}^{r-1}d_i^i}}\Big\}.\nonumber
\end{gather}
\end{satz*}

\begin{bew}{}
If $h$ is positive definite, we may use \ref{PosDefThm} with the values given in \ref{latticetopol} to obtain
$$\min_{h\in \LL^+}(f,h) \geq \frac{\sqrt{3}\sqrt[n]{\prod_{i=1}^{n-1} d_i^{n-i}}}{\dsum{1}{n-1}}\min_{h\in \LL_1}(f,h)
  = \frac{\sqrt{3}}{\sqrt[n]{\prod_{i=1}^{n-1} d_i^i}}\min_{h\in \LL_1}(f,h)$$
which proves \eqref{BnCpart1}.

The value $C(\LL)$ is constructed from terms that give valid bounds for the different possible cases $r:=\operatorname{rank}(h)=1,\dots,n$.
The first term, which is 1, obviously covers for $h$ of rank $r=1$. The term for $r=n$ has already been established in \eqref{BnCpart1}.

For positive semi-definite $h$ of rank $r$ with $1<r<n$, we proceed along the lines of Theorem~3 in \cite{BC}.

We can give $h$ as $h(x)=\transpose{x}Bx$ where $B$ is a rational singular matrix; the equation
$Bx=0$
hence has a rational solution $x\neq0$. Multiplying by a suitable rational number, we obtain a primitive integral vector $v=(v_1,\dots,v_n)$ with
$Bv=0.$
According to \ref{specialT} we can find an integral unimodular matrix $T$ of the form
$$T=\begin{pmatrix}
\ast&d_1&\dsum{1}{2}&\dots&\dsum{1}{n-2}&v_1\\
\ast&\ast&d_2&&\dsum{2}{n-2}&v_2\\
\vdots& &\ddots&\ddots&\vdots&\vdots\\
\ast& &\dots&\ast&d_{n-2}&v_{n-2}\\
\ast&&\dots&&\ast&v_{n-1}\\
\ast&&\dots&&\ast&v_n\\
\end{pmatrix}.$$
We now replace $f$ and $h$ by $\transpose{T}^{-1}f$ and $Th$, respectively; this leaves $M(f)$ and $(f,h)$ unchanged.
The matrix $B$ of $h$ is replaced by the matrix $\transpose{T}BT$ and, since $Bv=0$, the integral form $h$ has been replaced by an integral form in the $n-1$ variables $x_1,\dots,x_{n-1}$.
Furthermore, the special form of $T$ guarantees that $\transpose{T}BT\in \LL$.
We may clearly repeat this procedure until $h(x)$ is expressed as a positive
definite integral form in the variables $x_1,\dots,x_r$.
Let
\begin{align*}
\overline h(x_1,\dots,x_r) & := h(x)=h(x_1,\dots,x_r,0,\dots,0),\\
\overline f(x_1,\dots,x_r) & := f(x_1,\dots,x_r,0,\dots,0).
\end{align*}
Then $\overline f,\overline h$ are positive definite forms in $r$ variables, and $\overline h$ is integral. Clearly we have $M(\overline f)\geq M(f)$ and $(\overline f, \overline h)=(f,h).$
With respect to the sublattice
$$\overline \LL:=\left\{\begin{pmatrix}\Z&\dots&\dsum{1}{r-1}\Z&0&\dots&0\\\vdots&&\vdots&\vdots&&\\\dsum{1}{r-1}\Z&\dots&\dsum{1}{r-1}\Z&0&&\vdots\\0&\dots&0&0&&\\\vdots&&&&\ddots&\\0&&\dots&&&0\end{pmatrix}\right\}\cap \LL\subset \LL$$
(which contains $\overline h$) we may therefore use \eqref{BnCpart1} to obtain
$$\min_{\overline h\in\overline \LL\text{ of rank $r$}}(\overline f, \overline h)\geq\frac{\sqrt{3}}{\sqrt[r]{\prod_{i=1}^{r-1} d_i^i}}\min_{\overline h\in\overline \LL_1}(\overline f, \overline h).$$
Hence, we have
$$ (f,h) = (\overline f, \overline h)\geq\min_{\overline h\in\overline \LL^+}(\overline f,\overline h)
 \geq \frac{\sqrt{3}}{\sqrt[r]{\prod_{i=1}^{r-1} d_i^i}}\min_{\overline h\in\overline \LL_1}(\overline f, \overline h)
 \stackrel{\overline \LL_1\subset \LL_1}{\geq} \frac{\sqrt{3}}{\sqrt[r]{\prod_{i=1}^{r-1}d_i^{i}}}\min_{h\in \LL_1}(f,h).$$
This construction supplies all the other terms in $C(\LL)$ and thus ends the proof.
\end{bew}

\begin{anm}{}
\ref{BnCgeneral} can now be used as a substitute for \ref{BnCused}. This leads to the following generalisation of \ref{Tai1.1}:
\end{anm}
\begin{satz*}{}\label{Taigeneral}\\
Assume a (non-principal) polarisation $(1,d_1,\dots,\dsum{1}{g-1})$ and let $\LL$ be its Tits lattice. Suppose $\D=\S_g$ and let $\Gamma=\GampolConj$ or $\Gamma=\GampolConj(n)$ with $\gcd(n,\dsum{1}{g-1})=1$. Furthermore, let $\chi$ and $\omega$ be as in \ref{Tai1.1}. Then
$$\chi\omega^{\otimes l}\text{extends to }\overline{\D/\Gamma}^0
\iff
\left\{\begin{gathered}
\text{$\chi$ vanishes on all}\\\text{rational corank-1 boundary components}\\\text{of order at least $\nicefrac{l}{C(\LL)}$.}
\end{gathered}\right.$$
\end{satz*}

\section{How to get from $\A_g$ to $\Apol(n)$}
\label{genConstr}

Our main goal is to investigate the Kodaira dimension of $\Apol(n)$. Our method needs a non-trivial cusp form with respect to $\GampolConj(n)$ which we do not yet have.
However, the product $\chi$ of all even theta constants is a cusp form\footnote{see \cite[p.~42, Satz~3.3]{Freitag83}} with respect to $\Sp(2g,\Z)$. Denote its weight by $w_\chi$ and its order of vanishing on the cusp of $\A_g^\ast$ by $v_\chi$.
How can we use $\chi$ to construct a cusp form on $(\Apol(n))^\ast$?

\subsection{Maps, cusps and branching}

We have the following situation:
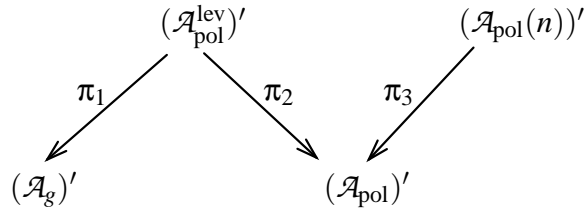
\begin{center}
\setlength{\unitlength}{0.00087489in}
\begingroup\makeatletter\ifx\SetFigFont\undefined%
\gdef\SetFigFont#1#2#3#4#5{%
  \reset@font\fontsize{#1}{#2pt}%
  \fontfamily{#3}\fontseries{#4}\fontshape{#5}%
  \selectfont}%
\fi\endgroup%
{\renewcommand{\dashlinestretch}{30}
\begin{picture}(2812,1219)(0,-10)
\thicklines
\path(2214.464,396.314)(2160.000,283.000)(2269.293,345.140)
\path(2160,283)(2790,958)
\path(1170,923)(1845,293)
\path(1731.686,347.464)(1845.000,293.000)(1782.860,402.293)
\path(945,913)(225,283)
\path(290.615,390.242)(225.000,283.000)(340.003,333.799)
\put(0,58){\makebox(0,0)[lb]{\smash{{{\SetFigFont{12}{14.4}{\rmdefault}{\mddefault}{\updefault}$(\A_g)'$}}}}}
\put(1890,58){\makebox(0,0)[lb]{\smash{{{\SetFigFont{12}{14.4}{\rmdefault}{\mddefault}{\updefault}$(\Apol)'$}}}}}
\put(900,1048){\makebox(0,0)[lb]{\smash{{{\SetFigFont{12}{14.4}{\rmdefault}{\mddefault}{\updefault}$(\Apollev)'$}}}}}
\put(2700,1048){\makebox(0,0)[lb]{\smash{{{\SetFigFont{12}{14.4}{\rmdefault}{\mddefault}{\updefault}$(\Apol(n))'$}}}}}
\put(405,643){\makebox(0,0)[lb]{\smash{{{\SetFigFont{12}{14.4}{\rmdefault}{\mddefault}{\updefault}$\pi_1$}}}}}
\put(1530,643){\makebox(0,0)[lb]{\smash{{{\SetFigFont{12}{14.4}{\rmdefault}{\mddefault}{\updefault}$\pi_2$}}}}}
\put(2250,643){\makebox(0,0)[lb]{\smash{{{\SetFigFont{12}{14.4}{\rmdefault}{\mddefault}{\updefault}$\pi_3$}}}}}
\end{picture}
}

\end{center}
where by $\A'$ we denote Mumford's partial compactification of $\A$. This is constructed from $\A$ by adding only the corank-1 boundary components. Note that this construction is well defined since it does not depend on a fan and that all these maps exist due to the inclusion relations of the corresponding groups.

What do we know about the partial compactifications of these spaces?
First of all, we know\footnote{See \cite[Part~I, Lemma~3.11]{HKW}} that $(\A_g)'$ has only a single cusp which we shall call $C_0$.

In $(\Apol)'$ there are several rational corank-1 boundary components which we shall denote by $C_1,\dots,C_u$.
Fix $i$ in $1,\dots,u$ and
denote the irreducible components of the reduction of $\pi_2^\ast C_i$ by $C_i^1,\dots,C_i^{v_i}\subset(\Apollev)'$. To each $C_i^j$ we can associate a unique\footnote{up to multiplication with $-1$} primitive vector in $\Z^{2g}$ that generates the corresponding isotropic space. By abuse of notation, we also denote this generator by $C_i^j$.
Let $\Cpollev(i)$ be a set of vectors that is a full system of representatives for these boundary components.
Denote the order of branching of $\pi_1$ in $C_i^j$ by $m_1(i,j)$
and that of $\pi_2$ by $m_2(i,j)$.

We know that $\GampollevConj$ is a normal subgroup of $\GampolConj$ and so
$\pi_2:\Apollev\rightarrow\Apol$ is a Galois cover.
The Galois group $\GamG:=\GampolConj/\GampollevConj$ operates transitively on $\Cpollev(i)$, so that for any fixed $i$
the order of the stabiliser $\Stab_{\GamG}(C):=\{g\in\GamG\where g(C)=C\}$ is the same for all $C\in\Cpollev(i)$.
If $-1\not\in\GamG$ (from \cite{Brasch93} we see that this is implied by $\dsum{1}{g-1}>2$), it can be given by
\begin{equation}\label{staborder}
|\Stab_{\GamG}(C_i^j)|=\frac{|\GamG|}{|\Cpollev(i)|}.
\end{equation}

Furthermore, the values $m_2(i,j)$ are the same for all $C\in\Cpollev(i)$ and we
can denote them by $m_2(i)$. So we have
$$\pi_2^\ast C_i = \sum_j m_2(i,j)C_i^j = m_2(i)\sum_j C_i^j.$$

\subsection{Modular forms}

From \refTits{GamSubSp} we know that $\GampollevConj\subset\Sp(2g,\Z)$ and hence $\chi$ is also a cusp form with respect to $\GampollevConj$.
On $C_i^j$ it vanishes of order $\ord(\chi,C_i^j)=v_\chi m_1(i,j)$.

Define $\chisym$ to be the symmetrisation of $\chi$ with respect to the Galois group $\GamG$ constructed as in \cite{H}.
This is a cusp form with respect to $\GampollevConj$ of weight $w_{\text{sym}}=|\GamG|w_\chi$.
We may choose any one cusp $C_i^1$ and
have $$\forall C\in\Cpollev(i)\suchthat\ord(\chisym,C)=\ord(\chisym,C_i^1).$$
To be precise, we have
\begin{align*}
\ord(\chisym,C_i^1)&=\sum_{a\in\GamG}\ord(\chi,a^{-1}(C_i^1))
 = \sum_{C_i^j\in\Cpollev(i)}|\Stab_{\GamG}(C_i^j)|\ord(\chi,C_i^j)\\
& = \sum_{C_i^j\in\Cpollev(i)}\frac{|\GamG|}{|\Cpollev(i)|}v_\chi m_1(i,j)
 = v_\chi\frac{|\GamG|}{|\Cpollev(i)|}\sum_{C_i^j\in\Cpollev(i)}m_1(i,j).
\end{align*}
For easier notation define $$M_1(i):=\sum_{C_i^j\in\Cpollev(i)}m_1(i,j).$$
In fact, $\chisym$ is also a cusp form with respect to $\GampolConj$. To make clear which group we are referring to we use the notation $\overline\chi$ in case of this second group.
On $(\Apol)'$ we now have
$$\ord(\overline\chi,C_i) = \ord(\chisym,C_i^1)/m_2(i)
 = v_\chi\frac{|\GamG|}{m_2(i)|\Cpollev(i)|}M_1(i).$$

\subsection{Vanishing on higher codimension}

So far we are able to control the order of vanishing on the corank-1 boundary components of a compactification of $\Apol$. This compactification may, however, be singular.
Assume we are given a $\GampolConj$-admissible collection of fans $\Sigma$ and obtain the corresponding compactification $(\Apol)^\ast$.
According to \cite[Theorem 7.20]{Nami} and \cite[Theorem 7.26]{Nami}, there exists a refinement $\tilde\Sigma$ of the collection $\Sigma$, which is also $\GampolConj$-admissible, such that the corresponding compactification $(\Apol)^\sim$ is stack-smooth. By this we mean that all fans are basic and hence no singularities arise from the toroidal construction but are only introduced by the group action.
Furthermore, we also get that the map $(\Apol)^\sim\rightarrow(\Apol)^\ast$ is a blowing-up and hence $(\Apol)^\sim$ is constructed from $(\Apol)^\ast$ by inserting new boundary divisors. These also correspond to rays in the closure $\overline C$ of the  cone of symmetric, positive definite matrices, as do the corank-1 boundary components, but here the rays are generated by matrices of rank strictly greater than 1.

We are now ready to proceed to the map $\pi_3$. Assume that the level $n$
is such that $\pi_3$ is branched of order $n$ along all boundary components.
For any cusp $C$ in the pullback $\pi_3^\ast C_i$ we then have
\begin{equation}\label{ordC}
\ord(\overline\chi,C) = n\ord(\overline\chi,C_i) = nv_\chi\frac{|\GamG|}{m_2(i)|\Cpollev(i)|}M_1(i).
\end{equation}
Now we use the generalised Barnes and Cohn \ref{BnCgeneral} on $(\Apol(n))^\sim$ which states that $\overline\chi$ vanishes on all of the boundary at least of order $\ord(\overline\chi,C)C(\LL).$ 
On the other hand, $\overline\chi$ is a modular form of weight $w_{\overline\chi}=|\GamG|w_\chi$ with respect to $\GampolConj(n)\subset\GampolConj$. This leads to the following equation for $(\Apol(n))^\sim:$
$$w_\chi|\GamG|L=\ord(\overline\chi,C)C(\LL)D+D_{\text{eff}}$$
where $L$ is the divisor corresponding to the ($\Q$-)line bundle\footnote{For $n\geq3$ this is in fact a line bundle.} of modular forms of weight 1 on $\Apol(n)$, $D$ is the boundary divisor of $(\Apol(n))^\sim$ and $D_{\text{eff}}$ is some effective divisor that we do not need to specify more precisely. This implies
$$-D=-\frac{w_\chi|\GamG|}{\ord(\overline\chi,C)C(\LL)}L+D'_{\text{eff}}.$$
Assume now that $n\geq3$ such that $\GampolConj(n)$ is neat. For any smooth toroidal compactification of $\Apol(n)$ we obtain
\begin{align*}
 K &= (g+1)L-D \\
  &=\Big[(g+1)-\frac{w_\chi|\GamG|}{\ord(\overline\chi,C)C(\LL)}\Big]L+D'_{\text{eff}}.
\end{align*}
We know from Mumford's extension of Hirzebruch proportionality (see \cite[Corollary~3.5]{Mu77}) that $h^0(L^k)\sim k^{\frac{1}{2}g(g+1)}$. We can therefore conclude that $h^0(K^k)\sim h^0(L^k)\sim k^{\frac{1}{2}g(g+1)}$ and hence that the Kodaira dimension is maximal if the coefficient of $L$ is positive.
This means we want
\begin{align}
 & \ord(\overline\chi,C)C(\LL)>\frac{w_\chi|\GamG|}{g+1}\nonumber\\
\iff & n\frac{v_\chi|\GamG|M_1(i)C(\LL)}{m_2(i)|\Cpollev(i)|} > \frac{w_\chi|\GamG|}{g+1}\nonumber\\
\iff & n>\frac{w_\chi m_2(i)|\Cpollev(i)|}{(g+1)v_\chi M_1(i)C(\LL)}.\label{ngthan3}
\end{align}

\section{Branching of the maps}

\subsection{Maps between toroidal varieties}

\begin{satz}{: Maps of toroidal varieties}
\label{ordBranch}
Assume we are given two arithmetic subgroups $\Gamma_1\subset\Gamma_2$ of $\Sp(2g,\Z)$ and a collection of fans $\tilde\Sigma$ that is admissible for both groups. Let $\A_i^\ast:=(\nicefrac{\S_g}{\Gamma_i})^\ast$. Then we have a map $\pi:\A_1^\ast\rightarrow\A_2^\ast$. Furthermore,
for a corank-1 boundary component $F$ the order of branching of $\pi$ on $F$ is given by the index $[\Pdash{\Gamma_2}{F}:\Pdash{\Gamma_1}{F}]$.
\end{satz}
\begin{bew}{}
The existence of $\pi$ follows easily from \cite[Theorem~1.13]{Oda} since
$$N'=\Pdash{\Gamma_1}{F}=\Pdash{\Sp(2g,\Z)}{F}\cap\Gamma_1\subset\Pdash{\Sp(2g,\Z)}{F}\cap\Gamma_2=\Pdash{\Gamma_2}{F}=N$$
where $[N:N']<\infty$ due to the choice of $\Gamma_i$.
We can glue the maps $\varphi_{F,\ast}$ since $\tilde\Sigma$ is admissible.

Since $F$ has corank 1, the groups $\Pdash{\Gamma_j}{F}$ for $j=1,2$ are 1-dimensional lattices. To ease the notation, we only consider the case $F=F_0$, but the construction goes through the same for all other rational corank 1 boundary components. For $F_0$, the quotient maps $e_j(F_0)$ are given by
$$e_j(F_0):\left\{\begin{array}{ccc}\S_g&\rightarrow&X_j(F_0)=\C^\ast\times\C^{g-1}\times\S_{g-1}\\
(\tau_{1,1},\tau_{1,2},\dots,\tau_{g,g})&\mapsto&(t_j,\tau_{1,2},\dots,\tau_{1,g},\tau')\end{array}\right.$$
where $\tau'=(\tau_{m,n})_{m,n\geq2}$ and $t_j=e^{2\pi i\tau_{1,1}/k_j}$ for some $k_j\in\N,j=1,2$.
Now we have a map
$$\tilde\pi:\left\{\begin{array}{ccl}X_1(F_0)&\rightarrow&X_2(F_0)\\
t_1&\mapsto&t_2=(t_1)^{\nicefrac{k_1}{k_2}}\\
\tau_{m,n}&\mapsto&\tau_{m,n}\quad\text{for all $(m,n)\neq(1,1)$}\end{array}\right..$$
This map extends naturally to the boundary $\{0\}\times\C^{g-1}\times\S_{g-1}$ of $\C^g\times\S_{g-1}$. Obviously, the order of branching of $\tilde\pi$ in $\{0\}\times\C^{g-1}\times\S_{g-1}$ is $\frac{k_1}{k_2}$.

Now we have to consider the quotient maps $q_j$
$$X_j(F_0)\hookrightarrow X_{\Sigma,j}(F_0)\stackrel{q_j}{\rightarrow}\nicefrac{X_{\Sigma,j}(F_0)}{P''_{\Gamma_j}(F_0)}\hookrightarrow\A_j^\ast.$$
According to \cite[Proposition~3.90 and Proposition~3.91]{HKW} the group $\cP''(F_0)$ can be identified as the group consisting of the block matrices
$$\begin{pmatrix}\epsilon&m&n\\0&A&B\\0&C&D\end{pmatrix}\in\GL(g+1,\R)$$
where $\smallsqmatrixtwo{A}{B}{C}{D}\in\Sp(2(g-1),\R)$, $\epsilon\in\R$ and $m,n\in\R^{g-1}$. The action of its generators
$$g''_1=\begin{pmatrix}1&0&0\\0&A&B\\0&C&D\end{pmatrix},\quad
g''_2=\begin{pmatrix}\epsilon&0&0\\0&\UnitMx&0\\0&0&\UnitMx\end{pmatrix}\quad\text{and}\quad
g''_3=\begin{pmatrix}1&m&n\\0&\UnitMx&0\\0&0&\UnitMx\end{pmatrix}$$
on $\tau=(\tau_1,\tau_2)\in\C^{g-1}\times\S_{g-1}$ is given by
\begin{align*}
g''_1(\tau)&=(\tau_1(C\tau_2+D)^{-1},(A\tau_2+B)(C\tau_2+D)^{-1})\\
g''_2(\tau)&=(\tau_1\epsilon,\tau_2)\\
g''_3(\tau)&=(\tau_1+m\tau_2+n,\tau_2).
\end{align*}
Now suppose that $g''=g''_1g''_2g''_3\in\cP''(F_0)$ is an element that operates like the identity on all of the boundary. Obviously, its action on the second component is determined by the submatrix $M:=\smallsqmatrixtwo{A}{B}{C}{D}$, and hence we need to have $M=\pm\UnitMx$. If $M=\UnitMx$, the factor $g''_1$ leaves $\tau_1$ invariant, and otherwise changes its sign. It is easy to see that in both cases $m=n=0$ and $\epsilon=\pm1$, where $\epsilon=1$ if and only if $M=\UnitMx$.
Hence, the only elements of $\cP''(F_0)$ that operate like the identity on all of the boundary are in fact $\UnitMx\in\cP''(F_0)$ and $-\UnitMx\in\cP''(F_0)$.
The same remains true if we intersect $\cP''(F_0)$ with the appropriate group $\Gamma_j$.
But since $-\UnitMx$ operates trivially on all of $\C^g\times\S_{g-1}$, this shows that the maps $q_j$ are not branched along the boundary divisor.

We obtain that the order of branching of $\pi:\A_1^\ast\rightarrow\A_2^\ast$ on the rational boundary components of corank~1 is also given by $\frac{k_1}{k_2}=[\Pdash{\Gamma_2}{F}:\Pdash{\Gamma_1}{F}]$.
\end{bew}

\subsection{The geometry of $\Apollev\rightarrow\A_g$ and $\Apollev\rightarrow\Apol$}

We shall now focus on the geometry of the maps $\pi_1$ and $\pi_2$.
In particular, we shall state a lemma on the order of branching
for these maps in each corank-1 boundary component of $(\Apollev)'$.

\begin{lem}{: Order of branching}
For a rational corank-1 boundary component $F\subset(\Apollev)'$ the orders of branching of the maps between the partial compactifications $\pi_1:(\Apollev)'\rightarrow(\A_g)'$ and $\pi_2:(\Apollev)'\rightarrow(\Apol)'$ are given by
\begin{align*}
m_1(C) & := [\Pdash{\Sp(2g,\Z)}{C}:\Pdash{\GampollevConj}{C}]\quad\mbox{and}\\
m_2(C) & := [\Pdash{\GampolConj}{C}:\Pdash{\GampollevConj}{C}],
\end{align*}
respectively, where $\Pdash{\Gamma}{C}:=\cP'(F)\cap\Gamma\subset\cP(F)$ is
the relevant lattice part of the stabiliser of $F$ with $C=V(F)$.
\end{lem}
\begin{bew}{}
This is a specialisation of \ref{ordBranch}.
\end{bew}

Let us now give the general outline of how we want to perform this calculation in both cases. We do the calculations that are the same for all cases over the rationals, and only then intersect with the four different groups.

The group $\Sp(2g,\Q)$ has only a single corank-1 boundary component, namely $C_0\widehat{=}(0,\dots,0,1)\in\Z^{2g}$,
and for this cusp \cite[Paragraph~3D]{HKW} shows that
$$\PdashQ:=\Pdash{\Sp(2g,\Q)}{C_0}=\Big\{\sqmatrixtwo{\UnitMx}{S}{0}{\UnitMx}\mbox{ where }S=\diag(0,\dots,0,s)\text{ and }s\in\Q\Big\}.$$
From this information we calculate the groups $\Pdash{\Gamma}{C}$ for the other $\Gamma\subset\Sp(2g,\Q)$ and any cusp $C$ as follows:
Since all cusps are conjugate with respect to $\Sp(2g,\Q)$, we can always find a matrix $M\in\Sp(2g,\Q)$ such that
\begin{equation}\label{defM}
C=C_0M.
\end{equation}
This implies
$$
\Pdash{\Gamma}{C} =\Pdash{\Sp(2g,\Q)}{C}\cap\Gamma
=\big(M^{-1}\PdashQ M\big)\cap\Gamma
$$
which leads to the following lemma:

\begin{lem*}{}
\begin{align*}
m_1(C) & = [M^{-1}\PdashQ M\cap\Sp(2g,\Z):M^{-1}\PdashQ M\cap\GampollevConj]\mbox{ and}\\
m_2(C) & = [M^{-1}\PdashQ M\cap\GampolConj:M^{-1}\PdashQ M\cap\GampollevConj].
\end{align*}
\end{lem*}

Note that the matrices $Q'\in\Pdash{\Sp(2g,\Q)}{C}$ have the form
$Q'=Q+\UnitMx$ where
\begin{equation}\label{defQ}
Q:=M^{-1}\sqmatrixtwo{0}{S}{0}{0}M.
\end{equation}
To intersect the group $\Pdash{\Sp(2g,\Q)}{C}$ with $\Gamma$ we only need to consider the conditions imposed on $Q$ by the appropriate lemma from section \refTits{CongSect}.

It is easy to see that the inverse of a matrix $M=\smallsqmatrixtwo{\alpha}{\beta}{\gamma}{\delta}\in\Sp(2g,\Q)$ is given by $M^{-1}=\smallsqmatrixtwo{\transpose{\delta}}{-\transpose{\beta}}{-\transpose{\gamma}}{\transpose{\alpha}}$ where $\alpha,\beta,\gamma,\delta\in\Q^{g\times g}$.
Split the vector representing the cusp into two vectors of length $g$ such that $C=(\mathfrak{c}_1,\mathfrak{c}_2)$. Then equation \eqref{defM} implies that
$\mathfrak{c}_1$ and $\mathfrak{c}_2$ are the last rows of the matrices $\gamma$ and $\delta$, respectively. Since the matrix $S$ has only one non-zero entry $s\in\Q$ we see that
\begin{equation}\label{Qprod}
Q=s\binom{\transpose{\mathfrak{\,c}_2}}{-\transpose{\mathfrak{\,c}_1}}(\mathfrak{c}_1,\mathfrak{c}_2).
\end{equation}

We shall now give the explicit calculation in the two cases separately.
Note that the classification of the cusps in \cite{Tits} is done with respect to the conjugate groups $\Gampol=R\GampolConj R^{-1}$ and $\Gampollev=R\GampollevConj R^{-1}$ where $R:=\diag(1,\dots,1,1,d_1,\dots,\dsum{1}{g-1}).$

\begin{lem*}{: Branching of $\pi_1$}\label{branchingpi1gen}\\
For a cusp $C_i^j=(\Dsum{1}{g-1},\Dsum{2}{g-1}a_2,\dots,a_g,0,\Dsum{2}{g-1}a_{g+2},\dots,a_{2g})\in\Apollev$ given with respect to $\Gampollev$ the
order of branching of $\pi_1:\Apollev\rightarrow\A_g$ is given by
$$m_1(C_i^j)=\gcd(\Dsum{1}{g-1},\Dsum{2}{g-1}a_2,\dots,D_{g-1}a_{g-1},a_g)^2.$$
\end{lem*}

\begin{bew}{}
Recall from \refTits{orbitsgamlev} that any cusp can be represented in the
form given in the statement. Since we want to work with $\GampollevConj$ rather than with $\Gampollev$ we have to multiply by $R$ and obtain
$$C_i^jR=(\Dsum{1}{g-1},\Dsum{2}{g-1}a_2,\dots,a_g,0,d_1\Dsum{2}{g-1}a_{g+2},\dots,\dsum{1}{g-1}a_{2g}).$$
In case this is not a primitive vector we divide by $k:=\gcd(\Dsum{1}{g-1},\dots,\dsum{1}{g-1}a_{2g})$ to obtain as representative $C$ of the cusp.

We define $Q$ as in \eqref{defQ} and can now proceed by asking when $Q+\UnitMx$ is in $\Pdash{\Gamma}{C}=\Pdash{\Sp(2g,\Q)}{C}\cap\Gamma$ for $\Gamma=\Sp(2g,\Z)$ or $\Gamma=\Gampollev$, respectively.

When taking the intersection of $\Pdash{\Sp(2g,\Q)}{C}$ with $\Sp(2g,\Z)$ the only condition is that the matrix $Q$ be integer.
The first entry of the $g+1$st row is given by
$q_{g+1,1}=-\Dsum{1}{g-1}^2k^{-2}s.$
Substitute $t:=-q_{g+1,1}\in\Z$. With this substitution,
the diagonal elements of the lower left quarter of $Q$ give rise to the necessary conditions
$t(\frac{a_i}{\Dsum{1}{i-1}})^2\in\Z$ for $i=2,\dots,g$.
Some straightforward calculation shows that these are also sufficient.
Hence,
\begin{equation}\label{tin1}
Q\in\Sp(2g,\Z) \quad\iff\quad t\in\Big(\frac{\Dsum{1}{g-1}}{\gcd(\Dsum{i}{g-1}a_i)_{i=1,\dots,g}}\Big)^2\Z.
\end{equation}

Since $\GampollevConj\subset\Sp(2g,\Z)$ we also get this condition for $\Pdash{\GampollevConj}{C}$ but in addition we have to consider \refTits{CongGampollevConj}.
The conditions of the upper right quarter of $Q$
lead to the necessary conditions $t\in(\tfrac{a_{g+1+i}}{\Dsum{1}{i}})^2\Z$
for all $i=1,\dots,g-1$.
Again, 
these imply all other conditions on $Q$ and hence lead to
$$Q\in\GampollevConj\quad\iff\quad t\in\Big(\frac{\Dsum{1}{g-1}}{\gcd(\Dsum{i}{g-1}a_i,\Dsum{i}{g-1}a_{g+i})_{i=1,\dots,g}}\Big)^2\Z\quad=\quad\Dsum{1}{g-1}^2\Z.$$
Combining this with equation \eqref{tin1} gives $m_1(C)=\gcd(\Dsum{i}{g-1}a_i)^2_{i=1,\dots,g}$ as claimed.
\end{bew}

\begin{lem*}{: Branching of $\pi_2$}\label{branchingpi2gen}\\
For a cusp $C_i^j=(\Dsum{1}{g-1},\Dsum{2}{g-1}a_2,\dots,a_g,0,\Dsum{2}{g-1}a_{g+1},\dots,a_{2g})\in\Apollev$ given with respect to $\Gampollev$ the
order of branching of $\pi_2:\Apollev\rightarrow\Apol$ is given by
$$m_2(C_i^j)=\Dsum{1}{g-1}.$$
\end{lem*}

\begin{bew}{}
Since $\GampollevConj$ is a normal subgroup of $\GampolConj$, the map $\pi_2$ induces
a Galois covering. This means that we may
restrict the investigation of the cusps $C_i^j\in(\Apollev)'$ for any $j$ to the primitive vector $C_i^0=(\Dsum{1}{g-1},\Dsum{2}{g-1}a_2,\dots,D_{g-1}a_{g-1},1,0,0,0)$.

As before we define $M$ and $Q$ by \eqref{defM} and \eqref{defQ}, respectively.
Again, we obtain conditions on $Q$ by intersecting $\Pdash{\Sp(2g,\Q)}{C}$ with $\Gamma$.

Since $\mathfrak{c}_2=0$, according to \eqref{Qprod} the only non-zero entries of $Q$ are in the lower left quarter.
\refTits{CongGampollevConj} states that for $Q+\UnitMx\in\GampollevConj$ these entries need to be integers. In particular, $q_{2g,g}=s\cdot 1\cdot 1=s\in\Z$. Since now $s\transpose{\mathfrak{c}_1}\mathfrak{c}_1$ is obviously an integer matrix we obtain the equivalence
$$Q+\UnitMx\in\GampollevConj\quad\iff\quad s\in\Z.$$

For $Q+\UnitMx\in\GampolConj$ we consider \refTits{CongGampolConj} where for the lower left quarter we find the condition
$$-s\transpose{\mathfrak{c}_1}\mathfrak{c}_1\in\PolMx^{-1}\DP=\begin{pmatrix}\Z&\Z&\Z&\dots&\Z\\\Z&\frac{1}{d_1}\Z&\frac{1}{d_1}\Z&&\frac{1}{d_1}\Z\\\Z&\frac{1}{d_1}\Z&\frac{1}{\dsum{1}{2}}\Z&&\frac{1}{\dsum{1}{2}}\Z\\\vdots&&&\ddots&\vdots\\\Z&\frac{1}{d_1}\Z&\frac{1}{\dsum{1}{2}}\Z&\dots&\frac{1}{\dsum{1}{g-1}}\Z\end{pmatrix}.$$
The condition on the top right matrix entry reads $q_{g+1,g}=s\Dsum{1}{g-1}\in\Z$ and hence we know that $s\in\frac{1}{\Dsum{1}{g-1}}\Z$ is a necessary condition.
Some straightforward calculation shows that it is in fact also sufficient.
Therefore,
$$m_2(C_j^0)=[\Pdash{\GampolConj}{C_j^0}:\Pdash{\GampollevConj}{C_j^0}]=[\tfrac{1}{\Dsum{1}{g-1}}\Z:\Z]=\Dsum{1}{g-1}$$
which completes the proof.
\end{bew}

\subsection{Branching of $\Apol(n)\rightarrow\Apol$}

\begin{lem}{}
\label{Branchingpi3}
Assume $\gcd(n,\dsum{1}{g-1})=1$. Then $\pi_3:\Apol(n)\rightarrow\Apol$ is branched of order $n$ on all corank-1 boundary components.
\end{lem}
\begin{bew}{}
Let $D$ be a corank-1 boundary divisor. Denote the stabilisers of the corresponding isotropic line in the groups $\Gampol$ and $\Gampol(n)$ by $\Stab_{\Gampol}(D)$ and $\Stab_{\Gampol(n)}(D)$, respectively.
Since $D$ has corank 1, these stabilisers are one-dimensional lattices and can therefore be given by $\Stab_{\Gampol}(D)\simeq k_1\Z$ and $\Stab_{\Gampol(n)}(D)\simeq k_2\Z$. Since $\Gampol(n)\subset\Gampol$ by definition, we know that $k_1|k_2$. Since $\gcd(n,\dsum{1}{g-1})=1$, the congruence condition imposed by $\Gampol(n)$ implies that $k_2/k_1=n$ for every such pair of lattices. But this index is exactly the order of branching, which proves the claim.
\end{bew}

\section{General Type results}

Before we can proof our main theorem, we only need two more things:
First, we need a way to restrict the scope to square-free polarisations, i.\,e.\ polarisations of type $(1,d_1,\dots,\dsum{1}{g-1})$ where all $d_i$ are square-free. The following lemma will make this precise.

\begin{lem}{: Square-free}
\label{squarefree}
Let $(1,e_1,\dots,\esum{1}{g-1})$ be the type of a polarisation
where $e_i=d_is_i^2$ and all $d_i$ are square-free. Then $(1,d_1,\dots,\dsum{1}{g-1})$ is the type of a square-free polarisation. 
Let $S:=\diag(1,s_1,\dots,\ssum{1}{g-1}), T:=\smallsqmatrixtwo{S}{0}{0}{S^{-1}}$ and $U:=\smallsqmatrixtwo{S}{0}{0}{S}$.
Then we have
$$T^{-1}\GpC{e}T\subset\GpC{d},\quad T^{-1}\GplC{e}T\subset\GplC{d}$$
and
$$U^{-1}\Gp{e}U\subset\Gp{d},\quad U^{-1}\Gpl{e}U\subset\Gpl{d}.$$
\end{lem}

\begin{bew}{}
We use the description of the groups given in \cite{Tits} and define $\De$ and $\Dd$ accordingly.

Let us begin with the relation $T^{-1}\GpC{e}T\subset\GpC{d}$.
Denote the matrices for the two polarisations by $\PolMx_e:=\diag(1,e_1,\dots,\esum{1}{g-1})\quad\text{and}\quad\PolMx_d:=\diag(1,d_1,\dots,\dsum{1}{g-1}).$ We have
$\PolMx_d=S^{-1}\PolMx_eS^{-1}$. Furthermore, for $M\in\De$ we obtain by simple computation that $S^{-1}MS\in\Dd$. Therefore,
\begin{equation}\label{SESinD}
S^{-1}\De S\subset\Dd.
\end{equation}
For $\GpC{e}$, \refTits{CongGampolConj} tells us
\begin{align*}
T^{-1}\GpC{e}T & \subset \sqmatrixtwo{S^{-1}}{}{}{S}\sqmatrixtwo{\De}{\De\PolMx_e}{\PolMx_e^{-1}\De}{\PolMx_e^{-1}\De\PolMx_e}\sqmatrixtwo{S}{}{}{S^{-1}} \\
 & = \sqmatrixtwo{S^{-1}\De S}{(S^{-1}\De S)(S^{-1}\PolMx_eS^{-1})}{(S\PolMx_e^{-1}S)(S^{-1}\De S)}{(S\PolMx_e^{-1}S)(S^{-1}\De S)(S^{-1}\PolMx_eS^{-1})}\\
 & \stackrel{\eqref{SESinD}}{\subset}\sqmatrixtwo{\Dd}{\Dd\PolMx_d}{\PolMx_d^{-1}\Dd}{\PolMx_d^{-1}\Dd\PolMx_d}.
\end{align*}
Since on the other hand $\GpC{e}\subset\Sp(2g,\Q)$ and $T,T^{-1}\in\Sp(2g,\Q)$, we may use \refTits{CongGampolConj} to conclude
$$T^{-1}\GpC{e}T\subset\sqmatrixtwo{\Dd}{\Dd\PolMx_d}{\PolMx_d^{-1}\Dd}{\PolMx_d^{-1}\Dd\PolMx_d}\cap\Sp(2g,\Q)=\GpC{d}.$$

For the relation $T^{-1}\GplC{e}T\subset\GplC{d}$ we first note that the first relation we proved implies that $$T^{-1}\GplC{e}T\subset T^{-1}\GpC{e}T\subset\GpC{d},$$
so that we only need to show the additional conditions imposed by \refTits{CongGampollevConj}.
This lemma states that the matrices $M\in T^{-1}\GplC{e}T$ are those matrices of $\GpC{e}$ that have the form
\begin{gather*}
M\in\sqmatrixtwo{S^{-1}}{}{}{S}\left(\binom{\transpose{\mathfrak{\,e}}}{1_g}(1_g,\mathfrak{e})\otimes\Z+\UnitMx\right)\sqmatrixtwo{S}{}{}{S^{-1}}
=\binom{S^{-1}\transpose{\mathfrak{\,e}}}{S}(S,\mathfrak{e}S^{-1})\otimes\Z+\UnitMx
\end{gather*}
where $\mathfrak{e}=(1,e_1,\dots,\esum{1}{g-1})$. Since $e_is_i^{-1}=d_is_i$ for all $i=1,\dots,g-1$ this means
\begin{gather*}
M\in\binom{S\transpose{\mathfrak{\,d}}}{S}(S,\mathfrak{d}S)\otimes\Z+\UnitMx
\subset\binom{\transpose{\mathfrak{\,d}}}{1_g}(1_g,\mathfrak{d})\otimes\Z+\UnitMx
\end{gather*}
with $\mathfrak{d}=(1,d_1,\dots,\dsum{1}{g-1}).$
Hence, all these matrices also satisfy the conditions of $\GplC{d}$.

The other two relations follow from these by conjugating with $R$ as in the previous section.
\end{bew}

Second, we need two number theoretic functions for counting the rational boundary components.
\begin{defi}{: Generalised phi function and Sigma functions}
\label{defphi}
Let $n,k\in\N$ and $\alpha\in\C$.
Define
$$\varphi_k(n):=\big|\{(x_1,\dots,x_k)\in\Z_n^k\where\gcd(x_1,\dots,x_k,n)=1\}\big|\quad\text{and}\quad\sigma_\alpha(n):=\sum_{d|n}d^\alpha.$$
The function $\varphi_1$ is known as the {\em Euler phi function} which we also denote by $\varphi$.
The function $\sigma_0$ is known as the function $\tau$ that gives the number of divisors.
\end{defi}

\begin{lem}{}
The functions $\varphi_k$ and $\sigma_\alpha$ are multiplicative.
\end{lem}
\begin{bew}{}
A generalisation of \cite[Theorem~2.7]{Nathanson} can be used to show this for $\varphi_k$. For $\sigma_\alpha$ this is some straightforward computation.
\end{bew}

We shall now combine the facts collected so far to prove our main theorem.

\begin{satz}{: General type for general genus}
\label{GeneralMainTheorem}
For any genus $g\geq3$ and coprime $d_1,\dots,d_{g-1}\in\N$ with $\dsum{1}{g-1}\neq2$, the moduli space $\Apol(n)$ of 
$(1,d_1,\dots,\dsum{1}{g-1})$-polarised Abelian varieties with a full level-$n$ structure is of general type, provided $\gcd(n,\dsum{1}{g-1})=1, n\geq3$ and
$$ n > \frac{(2^g+1)\dsum{2}{g-2}}{(g+1)2^{g-3}}\min\left\{\frac{d_1}{C(\LL(d_1,\dots,d_{g-1}))},\frac{d_{g-1}}{C(\LL(d_{g-1},\dots,d_1))}\right\}$$
where $$C(\LL(x_1,\dots,x_{g-1}))=\min\left\{1,\min_{2\leq r\leq g}\Big\{\tfrac{\sqrt{3}}{\sqrt[r]{\prod_{i=1}^{r-1}x_i^i}}\Big\}\right\}.$$
\end{satz}
\begin{bew}{}
First of all, if $\dsum{1}{g-1}=1$ we are in the principally polarised case and much weaker bounds than the one given are already known. Hence, we may assume $\dsum{1}{g-1}>2$.

Since we have $n\geq3$ we know that $\GampolConj(n)$ is neat and hence operates without fixed points. This implies that the quotient by $P''$ introduces no singularities, and since $(\Apol)^\sim$ is stack-smooth we know that $(\Apol(n))^\sim$ is smooth.

Furthermore, we may assume the $d_i$ to be square-free. Otherwise, we may write $d_i=s_i^2e_i$ where the $e_i$ are square-free. Then, according to \ref{squarefree}, we can conjugate $\GpC{d}(n)$ such that it becomes a subgroup of $\GpC{e}(n)$. This means that we have a rational map $\pi_4:(\A_{\text{pol},d}(n))^\sim\rightarrow(\A_{\text{pol},e}(n))^\sim$ and after some blowing-up this map becomes a morphism. By this morphism each form on $(\A_{\text{pol},e}(n))^\sim$ gives rise to a form on a suitable blow-up of $(\A_{\text{pol},d}(n))^\sim$ which implies that, if we can show general type for the (square-free) polarisation $e$, we also have general type for the polarisation $d$.

We consider the construction given in section \ref{genConstr}. For $\A_g$, we know from \cite[p.~42, Satz~3.3]{Freitag83} and \cite[Theorem~2.10]{MuOpenProb} that we have a cusp form $\chi$ of weight $w_\chi=(2^g+1)2^{g-2}$ that vanishes of order $v_\chi=2^{2g-5}$.

The map $\pi_3$ needs to be branched of order $n$. According to \ref{Branchingpi3} this is implied by the condition $\gcd(n,\dsum{1}{g-1})=1$.
We can now calculate a bound for the level $n$ by the construction described in section \ref{genConstr}, which gives
\begin{equation}\label{ngen}
n>\frac{w_\chi m_2(i)|\Cpollev(i)|}{(g+1)v_\chi M_1(i)C(\LL)}.
\end{equation}
(At this point we need the fact that the polarisation is coprime.)
Let us now calculate this value explicitely.

From \refTits{orbitsgam} we know that the cusps of $\Apol$ are given by vectors of the form $$C_i=(\Dsum{1}{g-1},\Dsum{2}{g-1},\dots,D_{g-1},1,0,\dots,0).$$
Let us consider such a cusp and the set $\Cpollev(i)$ consisting of the primitive vectors of the form $$C_i^j=(\Dsum{1}{g-1},\Dsum{2}{g-1}a_2,\dots,a_g,0,\Dsum{2}{g-1}a_{g+2},\dots,a_{2g})$$
with $0\leq a_k,a_{g+k}<\Dsum{1}{k-1}$ for $k=2,\dots,g$. From \ref{branchingpi1gen} and \ref{branchingpi2gen} we know that
$$m_1(i,j) =\gcd(\Dsum{1}{g-1},\Dsum{2}{g-1}a_2,\dots,a_g)^2\quad\text{and}\quad m_2(i) = \Dsum{1}{g-1}.$$
Define $B_k|d_k$ for $k=1,\dots,g-1$ by $\Bsum{1}{k}^2=m_1(i,j)$. This definition is unique because the $d_k$ are coprime. We now have
\begin{equation}\label{Bisgcd}
\Bsum{1}{g-1}=\gcd(\Dsum{1}{g-1},\Dsum{2}{g-1}a_2,\dots,a_g).
\end{equation}
We count these vectors using \ref{countGCD}: let both the $d_i$ and $c_i$ of the lemma be equal to $D_i$ and let the $b_i$ of the lemma be equal to $B_i$. Then we obtain that the number of $(g-1)$-tuples $(a_2,\dots,a_g)$ satisfying equation~\eqref{Bisgcd} is $\prod_{j=1}^{g-1}\varphi_{g-j}(\frac{D_j}{B_j})$.

On the other hand, $C_i^j$ is a primitive vector, so we have
\begin{align*}
1 &= \gcd(\Dsum{1}{g-1},\Dsum{2}{g-1}a_2,\dots,a_g,0,\Dsum{2}{g-1}a_{g+2},\dots,a_{2g}) \\
 &= \gcd(\Bsum{1}{g-1},\Dsum{2}{g-1}a_{g+2},\dots,a_{2g})
\end{align*}
and \ref{countGCD} (this time by letting also the $c_i$ of the lemma to be equal to $B_i$) states that we have a choice of $\prod_{j=1}^{g-1}\varphi_{g-j}(B_j)\big(\frac{D_j}{B_j})^{g-j}$ values for the $(g-1)$-tuple $(a_{g+2},\dots,a_{2g})$. So all in all
\begin{align*}
\Big|\left\{C_i^j\in\Cpollev(i):m_1(i,j)=\Bsum{1}{g-1}^2\right\}\Big|
 &= \prod_{j=1}^{g-1}\varphi_{g-j}(\tfrac{D_j}{B_j})\prod_{j=1}^{g-1}\varphi_{g-j}(B_j)\big(\tfrac{D_j}{B_j}\big)^{g-j} \\
 &= \prod_{j=1}^{g-1}\varphi_{g-j}(D_j)\big(\tfrac{D_j}{B_j}\big)^{g-j}
\end{align*}
where we use the property that the $d_k$ and hence the $D_k$ are square-free for the multiplicativity of the functions $\varphi_{g-j}$.
Taking the unweighted and weighted sum over all $B_k|D_k$ we therefore get
\allowdisplaybreaks
\begin{align*}
|\Cpollev(i)| &= \sum_{B_1|D_1}\dots\sum_{B_{g-1}|D_{g-1}}\prod_{j=1}^{g-1}\varphi_{g-j}(D_j)\big(\tfrac{D_j}{B_j}\big)^{g-j} \\
 & = \Big(\prod_{j=1}^{g-1}\varphi_{g-j}(D_j)\Big)\sum_{B_1|D_1}\dots\sum_{B_{g-1}|D_{g-1}}\prod_{j=1}^{g-1}\big(\tfrac{D_j}{B_j}\big)^{g-j} \\
 &= \Big(\prod_{j=1}^{g-1}\varphi_{g-j}(D_j)\Big)\prod_{j=1}^{g-1}\sum_{B_j|D_j}\big(\tfrac{D_j}{B_j}\big)^{g-j}
 = \prod_{j=1}^{g-1}\Big(\varphi_{g-j}(D_j)\sum_{B'_j|D_j}(B'_j)^{g-j}\Big) \\
 & = \prod_{j=1}^{g-1}\varphi_{g-j}(D_j)\sigma_{g-j}(D_j)
\intertext{and analogously}
M_1(i) &= \sum_{B_1|D_1}\dots\sum_{B_{g-1}|D_{g-1}}\prod_{j=1}^{g-1}\varphi_{g-j}(D_j)\big(\tfrac{D_j}{B_j}\big)^{g-j}\cdot\Bsum{1}{g-1}^2
 = \prod_{j=1}^{g-1}\varphi_{g-j}(D_j)\prod_{j=1}^{g-1}\sum_{B_j|D_j}\tfrac{D_j^{g-j}}{B_j^{g-j-2}} \\
 &= \prod_{j=1}^{g-1}\varphi_{g-j}(D_j)\Big[\Big(\prod_{j=1}^{g-2}\sum_{B_j|D_j}\big(\tfrac{D_j}{B_j}\big)^{g-j-2}D_j^2\Big)\sum_{B_{g-1}|D_{g-1}}D_{g-1}B_{g-1}\Big] \\
 &= \prod_{j=1}^{g-1}\varphi_{g-j}(D_j)\Big[\Dsum{1}{g-2}^2D_{g-1}\Big(\prod_{j=1}^{g-2}\sigma_{g-j-2}(D_j)\Big)\sigma_1(D_{g-1})\Big].
\end{align*}
Inserting this into condition \eqref{ngen} (using $m_2(i)=\Dsum{1}{g-1}$) the product of the $\varphi_{g-j}$ cancels and we are left with
\begin{align*}
n &> \frac{w_\chi}{(g+1)v_\chi C(\LL)}\frac{\Dsum{1}{g-1}\cdot\prod_{j=1}^{g-1}\sigma_{g-j}(D_j)}{\Dsum{1}{g-2}^2D_{g-1}\sigma_1(D_{g-1})\prod_{j=1}^{g-2}\sigma_{g-j-2}(D_j)} \\
 &= \frac{w_\chi}{(g+1)v_\chi C(\LL)}\frac{\prod_{j=1}^{g-1}\sigma_{g-j}(D_j)}{\Dsum{1}{g-2}\sigma_1(D_{g-1})\prod_{j=1}^{g-2}\sigma_{g-j-2}(D_j)}
\intertext{and since $\sigma_{a+b}(D)=\sum_{B|D}B^{a+b}\leq\sum_{B|D}B^aD^b=D^b\sigma_a(D)$ this is implied by}
\Longleftarrow\quad n &> \frac{w_\chi}{(g+1)v_\chi C(\LL)}\frac{\prod_{j=1}^{g-2}\sigma_{g-j-2}(D_j)D_j^2\cdot\sigma_{g-(g-1)}(D_{g-1})}{\Dsum{1}{g-2}\sigma_1(D_{g-1})\prod_{j=1}^{g-2}\sigma_{g-j-2}(D_j)} \\
 &= \frac{w_\chi}{(g+1)v_\chi C(\LL)}\Dsum{1}{g-2}.
\end{align*}
This condition has to hold true for all valid $D_k|d_k$ which obviously gives the condition
\begin{align*}
n&>\frac{w_\chi\dsum{1}{g-2}}{(g+1)v_\chi C(\LL)}
 = \frac{(2^g+1)2^{g-2}\dsum{1}{g-2}}{(g+1)2^{2g-5} C(\LL)}
 = \frac{(2^g+1)\dsum{1}{g-2}}{(g+1)2^{g-3} C(\LL)}\\
 &= \frac{(2^g+1)\dsum{1}{g-2}}{(g+1)2^{g-3}\min\{\min_{2\leq r\leq g}\sqrt{3}\sqrt[r]{\prod_{i=1}^{r-1}d_i^i},1\}}.
\end{align*}
Finally, we may use the symmetry given in \cite{BL} to obtain the other term of the statement.
\end{bew}

To conclude this paper, we give the bound for some special kinds of polarisations as corollaries:

\begin{kor}{}
For any genus $g\geq3$ and $d\in\N$, $d\geq 3$, the moduli space $\Apol(n)$ of $(1,\dots,1,d)$-polarised Abelian varieties with a full level-$n$ structure is of general type, provided $\gcd(n,d)=1$, $n\geq3$ and
$$n>\frac{2^g+1}{(g+1)2^{g-3}\sqrt{3}}\sqrt[g]{d^{g-1}}.$$
The same bound for the level applies for the moduli space of $(1,d,\dots,d)$-polarised abelian varieties with a full level-$n$ structure.

If the polarisation is of type $(1,\dots,1,d,\dots,d)$ where $1<i<g-1$ is the number of 1's, the bound is
$$n>\frac{2^g+1}{(g+1)2^{g-3}\sqrt{3}}d\min\{1,\sqrt[g]{d^{\min\{i,g-i\}}}\}.$$
\end{kor}
\begin{bew}{}
These statements follow easily from \ref{GeneralMainTheorem} by explicitely determining the minima.
\end{bew}

\begin{anm}{}
To make this result more accessible, we give a table for the lower bounds for $n$ in the case of polarisations of type $(1,\dots,1,d)$. Note that we have disregarded the condition $\gcd(d,n)=1$ to make make the pattern more obvious.
$$\begin{array}{c|cccccccccccccccccc}
g\backslash d&3&4&5&6&7&8&9&10&11&12&13&14&15&16&17&18&19&20\\
\hline
3 &3&4&4&5&5&6&6&7&7&7&8&8&8&9&9&9&10&10\\
4 &3&3&4&4&5&5&6&6&6&7&7&8&8&8&9&9&9&10\\
5 &3&3&3&4&4&5&5&6&6&6&7&7&7&8&8&9&9&9\\
6 &3&3&3&3&4&4&5&5&5&6&6&7&7&7&8&8&8&9\\
7 &3&3&3&3&4&4&4&5&5&5&6&6&6&7&7&7&8&8\\
8 &3&3&3&3&3&4&4&4&5&5&5&6&6&6&7&7&7&8\\
9 &3&3&3&3&3&3&4&4&4&5&5&5&6&6&6&7&7&7\\
\hline
\end{array}$$
\end{anm}

\begin{kor}{}
Let $s,t>1$ be integers with $\gcd(s,t)=1$. Then the moduli space $\Apol(n)$ of $(1,s,st)$-polarised Abelian varieties with a full level-$n$ structure is of general type provided $\gcd(n,st)=1$ and
$$n>\tfrac{3}{4}\sqrt{3}\sqrt[3]{s^2t^2\min\{s,t\}^2}.$$
\end{kor}
\begin{bew}{}
Again, calculation of the minima in \ref{GeneralMainTheorem} leads to this statement.
\end{bew}

\begin{anm}{}
To give an impression of the case $g=3$, we give the following table of minimal values of $n$ for a fixed polarisation of type $(1,s,st)$ for arbitrary $s,t\in\N$. Where the level had to be increased to satisfy the condition $\gcd(n,st)=1$ this is denoted by a pair of brackets around the increased value. The empty spaces result from the conditions $st\neq2$ and $\gcd(s,t)=1$.
$$\begin{array}{c|ccccccccccccccc}
s\backslash t & 1&2&3&4&5&6&7&8&9&10\\
\hline
1 & 3& &(4)&(5)& 4& 5& 5&(7)&(7)& 7\\
2 & & & 7& &(11)& &(13)& &(17)& \\
3 &(4)& 7& &(17)& 17& &(22)& 23& &(29)\\
4 &(5)& &(17)& &(27)& & 31& &(37)& \\
5 & 4&(11)& 17&(27)& & 37& 41&(47)& 49& \\
6 & 5& & & & 37& &(53)& & & \\
7 & 5&(13)&(22)& 31& 41&(53)& &(71)& 76& 81\\
8 &(7)& & 23& &(47)& &(71)& &(91)& \\
9 &(7)&(17)& &(37)& 49& & 76&(91)& & 113\\
10 & 7& &(29)& & & & 81& & 113& \\
\hline
\end{array}$$
\end{anm}

\begin{anm}{}
For $g\geq16$ the bound given by \ref{GeneralMainTheorem} is too high: Y.-S.~Tai showed in \cite{TaiEgloffstein} that for these cases no level structure is required, while we always obtain $n\geq3$.
We need this constant condition to know that $\Gamma(n)$ is neat and so we have no singularities coming from the group action. However, if we consider only principal polarisations, Tai showed in \cite{Tai} that for $g\geq5$ all singularities that occur on a suitable toroidal compactification are canonical. An argument by Salvetti Manni\footnote{given in \cite[p.~19]{HS}} shows that a similar reasoning can be applied to $g=4$. Since this means that we can extend pluricanonical forms to a smooth model we may drop the condition $n\geq3$ in this case. The same reasoning we employed for \ref{GeneralMainTheorem} now leads to the bound 
$$n>\frac{2^g+1}{(g+1)2^{g-3}}$$
which gives (including the known results for $g=1,2$)
$$\begin{array}{c|ccccccccc}g&1&2&3&4&5&6&7&8&9\\\hline n&7&4&3&2&2&2&2&1&1\end{array}.$$
These are exactly the numbers given in \cite[p.~17]{HS}, except for $g=7$ where $n=1$ is known to be sufficient.
Note that for $g=1,2$ the above formula remains true and even gives a sharp bound.
Note also that this gives the known result that $\A_g$ is of general type for $g\geq8$. This was originally proved by E.~Freitag \cite{Freitag83}, respectively D.~Mumford \cite{MuOpenProb} and is better by 1 that the result by Y.-S.~Tai \cite{Tai}.
\end{anm}

\section{Appendix: Technical lemmata}

\begin{lem}{}
\label{specialT}
Let $v=(v_1,\dots,v_n)\in\Z^n$. Then we can find an integer matrix $T$ of the form
$$T=\begin{pmatrix}
\ast&\bullet&\bullet&\dots&\bullet&v_1\\
\ast&\ast&\bullet&&\bullet&v_2\\
\vdots& &\ddots&\ddots&\vdots&\vdots\\
\ast& &\dots&\ast&\bullet&v_{n-2}\\
\ast&&\dots&&\ast&v_{n-1}\\
\ast&&\dots&&\ast&v_n\\
\end{pmatrix}$$
(where the $\bullet$ are arbitrary fixed integer values) such that
$\det(T)=\gcd(v_1,\dots,v_n).$
\end{lem}
\begin{bew}{}
We prove the claim by induction.

For $n=2$ we have the matrix $T=\smallsqmatrixtwo{t_{11}}{v_1}{t_{21}}{v_2}.$ We can choose $t_{11},t_{21}$ such that
$$\det(T)=t_{11}v_2-t_{21}v_1=\gcd(v_1,v_2)$$
which completes this case.

Let $n\in\N$ be arbitrary and assume the claim holds for $n-1$. Let $T^{(i)}$ and $T^{(i,j)}$ denote the submatrices of $T$ that consist of the columns 2 to $n-1$ with the $i$th or $i$th and $j$th rows removed.
Expansion of the determinant along the 1st column shows that
$$\det(T)=t_{11}\begin{vmatrix} & v_2\\T^{(1)}&\vdots\\ &v_n\end{vmatrix}-t_{21}\begin{vmatrix}&v_1\\&v_3\\T^{(2)}&\vdots\\&v_n\end{vmatrix}\pm\dots-(-1)^n t_{n1}\begin{vmatrix}&v_1\\T^{(n)}&\vdots\\&v_{n-1}\end{vmatrix}.$$
In particular, $t_{11},\dots,t_{n1}$ can be chosen such that the claim holds if
\begin{equation}\label{defF}
F:=\gcd\left(\begin{vmatrix} & v_2\\T^{(1)}&\vdots\\ &v_n\end{vmatrix},\begin{vmatrix}&v_1\\&v_3\\T^{(2)}&\vdots\\&v_n\end{vmatrix},\dots,\begin{vmatrix} & v_1\\T^{(n)}&\vdots\\&v_{n-1}\end{vmatrix}\right)\stackrel{?}{=}\gcd(v_1,\dots,v_n).
\end{equation}
We now prove this equality.
The first matrix is a $n-1\times n-1$ matrix of the special form needed for the induction. We may therefore assume that
\begin{equation}\label{specialTind1}
\begin{vmatrix}&v_2\\T^{(1)}&\vdots\\&v_n\end{vmatrix}=\gcd(v_2,\dots,v_n)=:f.
\end{equation}
Furthermore, expansion of this determinant along the last column gives
\begin{align}
v_2|T^{(1,2)}|\mp\dots+(-1)^nv_n|T^{(1,n)}|&=-(-1)^n f \nonumber\\
\iff\tfrac{v_2}{f}|T^{(1,2)}|\mp\dots+(-1)^n\tfrac{v_n}{f}|T^{(1,n)}| &= -(-1)^n 1\nonumber\\
\implies \gcd\big(|T^{(1,2)}|,\dots,|T^{(1,n)}|\big) &= 1. \label{specialTind2}
\end{align}
Now we can simplify $F$ by using expansion of the determinants along the last columns. Almost all terms of these expansions are multiples of $f$ because they contain one of $v_2,\dots,v_n$. Since according to \eqref{specialTind1} the first term in the gcd of \eqref{defF} is equal to $f$ these terms are not needed to determine the value of $F$. Hence, we are left with
\begin{align*}
F &= \gcd\big(f, v_1|T^{(1,2)}|,\dots,v_1|T^{(1,n)}|\big)\\
 &= \gcd\Big(f, v_1\gcd\big(|T^{(1,2)}|,\dots,|T^{(1,n)}|\big)\Big)
 \stackrel{\eqref{specialTind2}}{=} \gcd(f,v_1) = \gcd(v_1,\dots,v_n).
\end{align*}
So, the equation in \eqref{defF} is true and hence we can find $T$ as claimed.
\end{bew}

\begin{lem}{}
\label{countGCD}
Let $k\in\N$ and let $d_1,\dots,d_k\in\N$ be coprime integers. Chose integers $c_1,\dots,c_k$ and $b_1,\dots,b_k$ satisfying $b_i|c_i|d_i$ for all $i=1,\dots,k$.
Then
$$\Big|\Big\{(x_1,\dots,x_k)\where 0\leq x_i<\dsum{1}{i}, \gcd(x_i\csum{i+1}{k})_{i=0}^k=\bsum{1}{k}\Big\}\Big|=\prod_{i=1}^k\varphi_{k+1-i}\Big(\frac{c_i}{b_i}\Big)\Big(\frac{d_i}{c_i}\Big)^{k+1-i}$$
where we let $x_0=1$ to ease the notation of the gcd.
\end{lem}
\begin{bew}{}
Define $\dsum{i}{k}^{(j)}:=\dsum{i}{j-1}\dsum{j+1}{k}$. Since the $d_i$ are coprime we can rewrite the $x_i$ as
$$x_i\equiv\sum_{j=1}^iy_{i,j}\dsum{1}{i}^{(j)}\mod\dsum{1}{i}\quad\text{for $i=2,\dots,k$}$$
where we may chose $0\leq y_{i,j}<d_j$ (and let $y_{1,1}:=x_1$). This, according to the Chinese Remainder Theorem, makes the $y_{i,j}$ unique. 
Now the condition we have to consider is
\begin{align*}
\bsum{1}{k} &= \gcd\Big(\csum{1}{k},\csum{2}{k}y_{1,1},\csum{3}{k}\sum_{j=1}^2y_{2,j}\dsum{1}{2}^{(j)},\dots,c_k\sum_{j=1}^{k-1}y_{k-1,j}\dsum{1}{k-1}^{(j)},\sum_{j=1}^ky_{k,j}\dsum{1}{k}^{(j)}\Big) \\
 &= \gcd\Big(\csum{1}{k},\csum{2}{k}y_{1,1},\sum_{j=1}^2\csum{1}{k}^{(j)}y_{2,j},\dots,\sum_{j=1}^k\csum{1}{k}^{(j)}y_{k,j}\Big)
\intertext{since $c_i|d_i$ and the $d_i$ are coprime. Furthermore, since $b_i|c_i$ and the $c_i$ are coprime we obtain that $b_j|y_{i,j}$ for all $1\leq j\leq i$. Now the above condition is equivalent to}
\iff \quad 1 &= \gcd\Big(\tfrac{\csum{1}{k}}{\bsum{1}{k}},\tfrac{\csum{2}{k}}{\bsum{2}{k}}\tfrac{y_{1,1}}{b_1},\sum_{j=1}^2\tfrac{\csum{1}{k}^{(j)}}{\bsum{1}{k}^{(j)}}\tfrac{y_{2,j}}{b_j},\dots,\sum_{j=1}^k\tfrac{\csum{1}{k}^{(j)}}{\bsum{1}{k}^{(j)}}\tfrac{y_{k,j}}{b_j}\Big).
\intertext{Now let $\tilde y_{i,j}:=y_{i,j}/b_j$. Then the equality above is equivalent to}
\iff \quad 1 &= \gcd\big(\tfrac{c_1}{b_1},\tilde y_{1,1},\dots,\tilde y_{k,1}\big)\cdot\ldots\cdot\gcd\big(\tfrac{c_k}{b_k},\tilde y_{k,k}\big) \\
\iff \quad 1 &= \gcd\big(\tfrac{c_j}{b_j},\tilde y_{j,j},\dots,\tilde y_{k,j}\big)\quad\forall j=1,\dots,k.
\end{align*}
Since we have chosen $0\leq y_{i,j}<d_j$ we know $0\leq\tilde y_{i,j}<\frac{d_j}{b_j}$. In the restricted range $0\leq\tilde y_{i,j}<\frac{c_j}{b_j}$ the number of possible $(k-j+1)$-tuples $(\tilde y_{j,j},\dots,\tilde y_{k,j})$ satisfying the conditions is given by $\varphi_{k-j+1}(\frac{c_j}{b_j})$. Since $c_j|d_j$ we have exactly $(\frac{d_j}{c_j})^{k-j+1}$ copies of this range. This gives the value claimed.
\end{bew}

\end{document}